\documentclass[reqno]{amsart}

\input xy
\xyoption{all}
\usepackage{accents}
\usepackage{epsfig}
\usepackage{color}
\usepackage{amsthm}
\usepackage{amssymb}
\usepackage{amsmath}
\usepackage{amscd}
\usepackage{amsopn}
\usepackage{graphicx}

\usepackage{xspace}

\usepackage{url}
\usepackage{enumitem, hyperref}\hypersetup{colorlinks}

\usepackage{color} %\textcolor{red}{Test}

\definecolor{darkred}{rgb}{1,0,0} %can change the intensity in [0,1]
\definecolor{darkgreen}{rgb}{0,0.8,0}
\definecolor{darkblue}{rgb}{0,0,1}

\hypersetup{colorlinks,
linkcolor=darkblue,
filecolor=darkgreen,
urlcolor=darkred,
citecolor=darkgreen}

\usepackage[OT2,T1]{fontenc}
\usepackage[russian,USenglish]{babel}

\makeatletter
\def\reflb#1#2{\begingroup
    #2%
    \def\@currentlabel{#2}%
    \phantomsection\label{#1}\endgroup
}
\makeatother

\numberwithin{equation}{section}
\newtheorem {Theorem}{Theorem}
\numberwithin{Theorem}{section}

\newtheorem {Question}[Theorem]    {Question}
\newtheorem {Conjecture}[Theorem]    {Conjecture}

\theoremstyle{definition}

\theoremstyle{remark}

\newtheorem{Example}[Theorem]{Example}

\def    \eps    {\epsilon}

\newcommand{\CS}{{\mathcal S}}

\newcommand{\id}{{\mathit id}}
\newcommand{\pt}{{\mathit pt}}
\newcommand{\const}{{\mathit const}}

\def    \R      {{\mathbb R}}

\def    \Z      {{\mathbb Z}}
\def    \N      {{\mathbb N}}
\def    \Q      {{\mathbb Q}}

\def    \T      {{\mathbb T}}
\def    \CP     {{\mathbb C}{\mathbb P}}

\def    \12    {{\frac{1}{2}}}

\def    \p      {\partial}

\def    \H     {\operatorname{H}}

\def    \Fix     {\operatorname{Fix}}

\def    \width    {\operatorname{width}}

\def    \s  {\operatorname{c}}

\def \inv   {\mathrm{inv}}

\begin{document}

%%%%%%%%%%%%%%%%%%%%%%%%%%%%%%
%   TEXT FORMATTING

\setlength{\smallskipamount}{6pt}
\setlength{\medskipamount}{10pt}
\setlength{\bigskipamount}{16pt}

%%%%%%%%%%%%%%%%%%%%%%%%%%

%%%%%%%%%%%%%%%%%%%%%%%%%%

%%%%%%%%%%%           BEGINNING OF  TEXT

%%%%%%%%%%%%%%%%%%%%%%%%%%

\title[Approximate Identities and Lagrangian Poincar\'e
Recurrence]{Approximate Identities and Lagrangian Poincar\'e Recurrence}

\author[Viktor Ginzburg]{Viktor L. Ginzburg}
\author[Ba\c sak G\"urel]{Ba\c sak Z. G\"urel}

\address{BG: Department of Mathematics, University of Central Florida,
  Orlando, FL 32816, USA} \email{basak.gurel@ucf.edu}

\address{VG: Department of Mathematics, UC Santa Cruz, Santa Cruz, CA
  95064, USA} \email{ginzburg@ucsc.edu}



\date{\today} 

\thanks{The work is partially supported by NSF CAREER award
  DMS-1454342 and NSF grant DMS-1440140 through MSRI (BG) and Simons
  Collaboration Grant 581382 (VG)}

%\bigskip

\maketitle

\begin{center}
  
% % %%%%  remove for the arxiv version  %%%%
% \emph{\foreignlanguage{russian}{Достаньте двойные листочки}\\ 

%   To Rafail Kalmanovich Gordin on the occasion of his 70th birthday}
% % %%%%

  \emph{To Rafail Kalmanovich Gordin on the occasion of his 70th
    birthday}

\end{center}

\begin{abstract} In this note we discuss three interconnected problems
  about dynamics of Hamiltonian or, more generally, just smooth
  diffeomorphisms. The first two concern the existence and properties
  of the maps whose iterations approximate the identity map with
  respect to some norm, e.g., $C^1$- or $C^0$-norm for general
  diffeomorphisms and the $\gamma$-norm in the Hamiltonian case, and
  the third problem is the Lagrangian Poincar\'e recurrence
  conjecture.

\end{abstract}

%\maketitle

%\bigskip

\tableofcontents

%\newpage

\section{Introduction}
\label{sec:intro}
In this note we focus on three interconnected problems concerning the
dynamics of general smooth and, more specifically, Hamiltonian
diffeomorphisms. (A Hamiltonian diffeomorphism is the time-one map of
the Hamiltonian isotopy generated by a time-dependent Hamiltonian.)
The first two problems are about the existence and properties of the
maps whose iterations approximate the identity map with respect to
some norm, e.g., $C^1$- or $C^0$-norm for general diffeomorphisms and
the $\gamma$-norm in the Hamiltonian case. (The $\gamma$-norm is a
natural norm on the group of Hamiltonian diffeomorphisms, coming from
min-max critical values in Floer theory a.k.a.\ spectral invariants.)
In dynamical systems theory such maps are often called rigid, but we
prefer the more intuitive term \emph{approximate identity}, borrowed
from analysis. The third problem is the Lagrangian Poincar\'e
recurrence conjecture.

The notion of an approximate identity and several variants of the
definition are spelled out in Section \ref{sec:AI}. Approximate
identities can have interesting and non-trivial dynamics, and the main
examples of such maps are the so-called pseudo-rotations; see
\cite{A-Z, Br:Invent, FKr, GG:PR, GG:PRvR}. Yet, in some way, these
maps resemble actions of compact abelian groups on manifolds. For
instance, one may expect the fixed point set of an approximate
identity to be nowhere dense and the map near isolated fixed points to
satisfy some non-degeneracy conditions. This is the essence of the
first problem (Question \ref{quest:interior}).

The second problem (Question \ref{quest:sa-ai}) concerns the existence
of Hamiltonian diffeomorphisms which are $\gamma$-approximate
identities.  This question appears to be related to the Conley
conjecture asserting that for many closed symplectic manifolds every
Hamiltonian diffeomorphism has infinitely many un-iterated periodic
orbits. The conjecture is known to hold for a broad class of manifolds
(see \cite{GG:survey, GG:Rev}) including all symplectically aspherical
closed manifolds, but there are some exceptions. For instance, an
irrational rotation of $S^2$ about the $z$-axis has exactly two
periodic points -- the poles.  Although there is no established formal
connection between the Conley conjecture and approximate identities,
all known ``counterexamples'' to the Conley conjecture are
$\gamma$-approximate identities.  In particular, as proved in
\cite[Thm.\ 5.1]{GG:PR}, every pseudo-rotation of $\CP^n$, i.e., a
Hamiltonian diffeomorphism with exactly $n+1$ periodic points is a
$\gamma$-approximate identity. In the spirit of the Conley conjecture,
it is reasonable to expect that a symplectically aspherical manifold
does not admit compactly supported $\gamma$-approximate identities.

Finally, in Section \ref{sec:LPR} we turn to the Lagrangian Poincar\'e
recurrence conjecture according to which the images of a closed
Lagrangian submanifold $L$ under the iterates of a Hamiltonian
diffeomorphism $\varphi$ cannot be all disjoint from $L$, i.e.,
$L\cap \varphi^{k_i}(L)\neq\emptyset$ for some sequence
$k_i\to\infty$; see Conjecture \ref{conj:LPR}. Very little is known
about this conjecture and, at the time of writing, the only
established non-trivial case is when $\varphi$ is a pseudo-rotation of
$\CP^n$, \cite[Thm.\ 4.2]{GG:PR}.

Naturally, in such a short note we cannot possibly spell out all 
necessary definitions and background results. Section \ref{sec:AI} can
be, perhaps, accessible to a graduate student with sufficient
background in manifolds, group actions and basic dynamical systems
theory. However, in Sections \ref{sec:AI-Ham} and \ref{sec:LPR} we
occasionally make use of fairly advanced notions from symplectic
topology (Floer theory) and Hamiltonian dynamics. For these notions
and results we refer the reader to, e.g., \cite{HZ, MS, Sa}.

\medskip

\noindent{\bf Acknowledgements.} The authors are grateful to Mita
Banik, Bassam Fayad, Fr\'ed\'eric Le Roux, Leonid Polterovich, Sobhan
Seyfaddini, Egor Shelukhin and the referee for useful comments. Parts
of this work were carried out during the first author's visit to NCTS
(Taipei, Taiwan) and while the second author was in residence at MSRI,
Berkeley, CA, during the Fall 2018 semester.  The authors would like
to thank these institutes for their warm hospitality and support.

\section{Approximate identities and almost periodic maps}
\label{sec:AI}
Consider a class of compactly supported $C^\infty$-diffeomorphisms
$\varphi$ of a smooth manifold $M$ (e.g., all such diffeomorphisms or
compactly supported Hamiltonian diffeomorphisms when $M$ is
symplectic, etc.), equipped with some metric, e.g., the $C^0$- or
$C^1$- or $C^r$-metric or the $\gamma$-metric in the Hamiltonian case
-- see Section \ref{sec:AI-Ham}. The norm $\|\varphi\|$ is by
definition the distance from $\varphi$ to the identity. A map
$\varphi$ is said to be a $\|\cdot\|$-\emph{approximate identity}, or
a $\|\cdot\|$-\emph{a.i.}\ for the sake of brevity, if
$\varphi^{k_i}\to \id$ with respect to the norm $\|\cdot\|$ for some
sequence $k_i\to\infty$. (Strictly speaking, the entire sequence of
iterates $\varphi^{k_i}$ should be called an approximate identity. We
believe that a confusion with approximate identities in analysis is
unlikely.) Clearly, a $C^1$-a.i.\ is automatically a $C^0$-a.i. The
definition extends to flows in an obvious way. Approximate identities
have been extensively studied in dynamics, although usually from a
perspective different than ours; see, e.g., \cite{A-Z, Br:Invent, FK,
  FKr, Ku} and references therein. In this section we focus on $C^0$-
and $C^1$-a.i.'s, but first some terminological remarks are due.

In dynamics, approximate identities are often called \emph{rigid}
maps. We find this terminology misleading, for the term ``rigid'' is
routinely used in a different sense. Furthermore, rigidity is often
associated with structural stability and we are not aware of any
situation where an a.i.\ would be structurally stable. In fact, in
many instances it is not hard to show that an a.i.\ cannot be
structurally stable in a suitable class of maps. (Regarding
terminology, we also note that in low-dimensional dynamics $C^0$-rigid
maps are sometimes called recurrent, \cite{KP}.)

One should keep in mind that the above definition allows some room for
pathological behavior. For example, hypothetically it is possible that
$\|\varphi^{k'_i}\|\to\infty$ for some sequence $k'_i\to\infty$ for an
a.i.\ $\varphi$.  The definition can be and has been refined and
amended in several ways. For instance, one can impose conditions on
the rate of convergence of $\varphi^{k_i}\to\id$ (e.g., Diophantine
vs.\ Liouville) or on the arithmetic properties of the sequence
$$
K_\eps=\{k \mid \|\varphi^k\|<\eps\},
$$
e.g., that $K_\eps$ has positive density or contains infinitely many
primes, etc., significantly restricting the class of a.i.'s and their
possible dynamics; cf.\ \cite{FKr}.

One such refinement is of particular relevance to us. Namely,
$\varphi$ is called \emph{$\|\cdot\|$-almost periodic} if for every
$\eps>0$ the sequence $K_\eps$, where we now take $K_\eps\subset \Z$,
is \emph{quasi-arithmetic}, i.e., the difference between any two
consecutive terms is bounded by a constant (possibly depending on
$\eps$). Clearly, an almost periodic map is an a.i. Topological
dynamics of $C^0$-almost periodic maps is studied in detail in
\cite{GH}; see also \cite{BR}. Almost periodic maps are closely
related to compact group actions on $M$: $\varphi$ is $C^0$-almost
periodic if and only if the family $\{\varphi^k\}$ is equicontinuous
and thus generates a compact abelian group $G$ of (compactly
supported) homeomorphisms,~\cite{GH}.
  
\begin{Example}[Actions of compact Lie groups]
  \label{ex:torus}
  A translation, $x\mapsto x+\alpha$ where $\alpha\in G$, in the
  circle or torus $G$ is $C^1$-almost periodic. More generally,
  $\varphi$ is $C^1$-almost periodic whenever it topologically
  generates a $C^1$-action $G\times M\to M$ of a compact abelian Lie
  group $G$ on $M$, i.e., the subgroup $\{\varphi^k\mid k\in\Z\}$ is
  relatively compact in the group of $C^1$-diffeomorphisms, or,
  equivalently, $\varphi$ is an isometry with respect to some metric.
\end{Example}

In fact, the converse is also true, although this is ultimately a
non-trivial result closely related to the Hilbert--Smith conjecture
asserting that a locally compact group acting effectively on a
manifold is a Lie group; see, e.g., \cite{MZ,Pa,Ta} and references
therein. Namely, it is not hard to show that, when $\varphi$ is
$C^1$-almost periodic, the family $\{\varphi^k\}$ is precompact in the
$C^1$-topology and thus generates an action of a compact abelian group
$G$ by compactly supported $C^1$-diffeomorphisms of $M$. Then $G$ is
necessarily a Lie group (this is a deep result, \cite{RS}) and the
action map $G\times M\to M$ is $C^1$; see \cite{MZ} and references
therein. As a consequence, $\varphi$ is exactly as in Example
\ref{ex:torus}, and the dynamics of $C^1$-almost periodic maps is
rather boring.

A $C^0$-a.i.\ can be thought of as a map with simultaneous return
times: for every $\eps>0$ there exist infinitely many $k\in\N$ such
that $d(x,\varphi^k(x))<\eps$ for all $x\in M$, where $d$ is the
distance with respect to some fixed metric on $M$. (The key point here
is that $k$ is independent of $x$.)  In particular, $\varphi$ is not a
$C^0$-a.i.\ whenever there exists a non-recurrent point $x\in M$. This
observation provides numerous examples of maps which are not
$C^0$-a.i.'s. For instance, $\varphi$ is not a $C^0$-a.i.\ when it has
a non-elliptic periodic orbit. However, it is not clear to the authors
if, say, a volume-preserving $C^1$-a.i.\ must have zero topological
entropy and vanishing all Lyapunov exponents, although this seems
rather probable; see \cite{KH} for the definitions and also \cite{ACW}
for some related results. (For a $C^0$-almost periodic diffeomorphism,
vanishing of the topological entropy is easy to prove.) Furthermore,
the growth rate of $D\varphi^k$ is closely related to a.i.-type
features of $\varphi$; see \cite{FP, Po:Invent, Po:slow, PS} and
references therein.  Finally, note that obviously a $C^0$-a.i.\ cannot
be mixing or topologically mixing, and a.i.'s are often studied
in connection with mixing properties, \cite{A-Z, FKr}. For instance, the
horocycle flow is mixing and hence not a $C^0$-a.i.; we refer the
reader to \cite{Para} for the proof and \cite{Ma} for refinements and
generalizations of this result.

There are examples of $C^0$-a.i.'s with interesting dynamics. 

\begin{Example}[Pseudo-rotations in 2D]
  \label{Ex:PR-2D}
  Let $\varphi$ be a pseudo-rotation of the sphere $M=S^2$ or the
  closed disk $D^2$, i.e., an area-preserving diffeomorphism with
  exactly two periodic points for $S^2$ or one periodic point for
  $D^2$. (In both cases, all periodic points are necessarily fixed
  points.)  Already in this case, $\varphi$ can have rather surprising
  dynamics. For instance, there exist pseudo-rotations $\varphi$ with
  exactly three invariant ergodic measures. For $S^2$ these are the
  two fixed points and the area form, and for $D^2$ one has to replace
  one of the fixed points by the Lebesgue measure on $\p D^2$. The
  linearizations $D\varphi$ at the fixed points (and also
  $\varphi|_{\p D^2}$) are rotations in the same angle
  $\theta\not\in\pi\Q$, and most of the orbits of $\varphi$ are
  dense. Such maps are constructed in \cite{AK}; see also
  \cite{FK}. An analytic pseudo-rotation of $D^2$ is locally, (i.e.,
  near the fixed point) a $C^\infty$-a.i., \cite{A-Z}. Furthermore,
  assume that the rotation number $\theta/\pi$ is \emph{exponentially
    Liouville}, i.e., for any $c>0$ it can be approximated by rational
  numbers $p/q$ with error smaller than $\exp(-c q)$. Then $\varphi$
  is also (globally) a $C^0$-a.i., \cite{Br:Invent, GG:PR}. At the
  same time, a pseudo-rotation with a dense orbit is never
  $C^0$-almost periodic. Indeed, focusing on $M=S^2$, assume the
  contrary and let $G$ be the compact topological abelian group
  generated by $\varphi$. Then the action of $G$ on $M$ is transitive
  since $\varphi$ has a dense orbit. Thus, by the corollary in
  \cite[Sect.\ 6.3, p.\ 243]{MZ}, $G$ is a compact abelian Lie
  group. However, it is clear that such a group cannot act
  transitively on~$S^2$.
\end{Example}

In some sense, approximate identities resemble actions of compact
abelian groups on manifolds, although Example \ref{Ex:PR-2D} shows
that an approximate identity can have intricate dynamics and need not
literally generate such an action.  In any event, our first question
is inspired by this analogy.  Let $\varphi\colon M\to M$ be a smooth,
compactly supported a.i.\ (in whatever sense) and let
$\Fix(\varphi)=\{x\in M\mid \varphi(x)=x\}$ be the set of its fixed
points.

\begin{Question}
  \label{quest:interior}
  Assume that $M$ is connected and $\varphi\neq \id$. Must
  $\Fix(\varphi)$ have empty interior? What can be said about
  $D_x\varphi\colon T_xM\to T_xM$ at $x\in \Fix(\varphi)$? For
  instance, can $D_x\varphi$ be degenerate when $x$ is an isolated
  point in $\Fix(\varphi)$?
\end{Question}

When $M$ is non-compact, the interior of $\Fix(\varphi)$ is necessarily
non-empty since $\varphi$ is compactly supported. Thus the affirmative
answer to the first part of the question would simply mean that such a
map does not exist in this case. It is an easy exercise to check this
for $n=1$. If $\varphi$ is $C^1$-almost periodic and thus generates
a compact Lie group $G$ of smooth transformations, the answers readily
follow from the discussion above and the standard fact that a
$G$-action can be linearized near a fixed point; see, e.g.,
\cite{GGK}.  When $\varphi$ is $C^0$-almost periodic, the negative
answer to the first question would follow from the Hilbert--Smith
conjecture and a theorem of Newman asserting that the interior of the
fixed point set $M^\Gamma$ is empty whenever $\Gamma$ is a finite
group acting on $M$ by homeomorphisms; see \cite{Ne} and also
\cite{Dr}. (The difficulty is that although $\varphi$ is $C^1$, we
have no control over the ``smoothness'' of the $G$-action, e.g., we do
not know that $G$ acts by Lipschitz transformations and thus is a Lie
group, cf.\ \cite{RS}.)  When $M$ is a closed surface (and in some
other instances in dimension two) and $\varphi$ is a $C^0$-a.i.\
homeomorphism homotopic to the identity, the interior of
$\Fix(\varphi)$ is empty; this follows from the results in
\cite{KP}. The question appears to be completely open in general, even
for $C^1$-a.i.'s and even when $M=\R^n$, $n\geq 3$.  Finally, we point
out that one can expect a.i.'s (in any sense) to be very non-generic.

\section{Approximate identities in the Hamiltonian setting}
\label{sec:AI-Ham}

In this section, which requires some background in symplectic
topology, we turn to the case where $(M^{2n},\omega)$ is a symplectic
manifold and $\varphi=\varphi_H$ is a Hamiltonian diffeomorphism of
$M$, generated by a time-dependent Hamiltonian
$H\colon S^1\times M\to \R$. For the sake of simplicity, the manifold
$M$ is assumed to be (positive) monotone or symplectically aspherical,
i.e., $\omega|_{\pi_2(M)}=0=c_1(TM)|_{\pi_2(M)}$, throughout the rest
of the paper. (We refer the reader to, e.g., \cite{HZ, MS, Sa} for the
necessary definitions.)  In this setting, there are two other natural
norms to consider in addition to the $C^1$- and $C^0$-norm. These are
the $\gamma$-norm and the Hofer norm. We will focus on the
former. When $M$ is closed,
$$
\gamma(\varphi)=\s_{[M]}(H)+\s_{[M]}\big(H^{\inv}\big),
$$
where $H^{\inv}$ is the Hamiltonian generating the flow
$\big(\varphi_H^{t}\big)^{-1}$ and $\s_w$ is the spectral invariant
associated with a class $w\in \H_*(M)$; see, e.g., \cite{EP, Oh:gamma,
  Oh:constr, Sc, Vi}. (For our purposes, it is convenient to think of
$\varphi_H$ as an element of the universal covering of the group of
Hamiltonian diffeomorphisms.) When $M$ is symplectically aspherical
and also for $M=\CP^n$ equipped hereafter with the standard symplectic
form,
$$\gamma(\varphi)=\s_{[M]}(H)-\s_{[\pt]}(H).
$$
One can also extend this definition to the case where $M$ is open,
provided that its structure at infinity is well-controlled, e.g., $M$
is convex at infinity; cf.\ \cite{FS, Gu07}. When $M$ is
symplectically aspherical, the $\gamma$-norm is continuous with
respect to the $C^0$-norm and thus a Hamiltonian $C^0$-a.i.\ is
automatically a $\gamma$-a.i, \cite{BHS}. This is also
true for $\CP^n$, \cite[Thm.\ C]{Sh}.

The authors learned of the following question from L. Polterovich and
Seyfaddini:

\begin{Question}
  \label{quest:sa-ai}
  Does a symplectically aspherical manifold $M$ admit compactly
  supported $C^0$- or $\gamma$-almost periodic Hamiltonian
  diffeomorphisms $\varphi\neq \id$ or, more generally, a.i.'s?
\end{Question}

It is reasonable to conjecture that the answer is negative in all
instances. As a warm-up, it is easy to see that $M$ does not admit a
$C^1$-a.i. In fact, $\|\varphi^k\|_{C^1}\to\infty$ even when
$\varphi\neq \id$ is degenerate.  Indeed, let $\CS(H)$ be the action
spectrum of $H$. This is a compact subset of $\R$, \cite{HZ}. Set
$\width(\varphi):=\max \CS(H)-\min\CS(H)$. Then $\width(\varphi)>0$
whenever $\varphi\neq \id$, \cite{Sc}. Hence,
$$
\width(\varphi^k)\geq k\width(\varphi)\to\infty.
$$
On the other hand, $\width(\varphi)\leq \const\|\varphi\|_{C^1}$.
Hence, $\width(\varphi^k)\leq O\big(\|\varphi^k\|_{C^1}\big)$.  (This
argument is taken from \cite{Po:Invent}.)

Furthermore, if we knew that a $C^0$-a.i.\ is automatically
non-degenerate (see Question \ref{quest:interior}), it would follow
that any closed symplectic manifold $M$ such that
$\H_{\mathit{odd}}(M;\Z)\neq 0$, e.g., a symplectically aspherical
manifold, does not admit Hamiltonian $C^0$-a.i.'s. (Indeed, by Floer
theory, every non-degenerate Hamiltonian diffeomorphism of $M$ must
have non-elliptic fixed points and thus cannot be a $C^0$-a.i.) Note
also that the closure of the group of Hamiltonian diffeomorphisms with
respect to the $\gamma$-norm is rather complicated and poorly
understood, but it certainly contains elements other than
homeomorphisms and does not act on $M$ in any obvious
sense,~\cite{Hu}.

One can expect few manifolds to admit $\gamma$-a.i.'s. Of
course, $M$ admits $C^1$-almost periodic Hamiltonian diffeomorphisms
(cf.\ Example \ref{ex:torus}) when it has a Hamiltonian
$S^1$-action. In particular, all symplectic toric manifolds and
coadjoint orbits of compact Lie groups (e.g., $\CP^n$, complex
Grassmannians, flag manifolds, etc.) admit $C^1$-almost periodic
Hamiltonian diffeomorphisms. However, there are no other known
manifolds having $\gamma$-a.i.'s. Overall, the situation seems to be
parallel to the Conley conjecture asserting that for many manifolds
(but obviously not all, e.g., $S^2$) every Hamiltonian diffeomorphism
has infinitely many un-iterated periodic orbits.  The conjecture has
been proved for a broad class of manifolds including all
symplectically aspherical manifolds; see \cite{GG:survey, GG:Rev}. All
known Hamiltonian diffeomorphisms with finitely many periodic orbits
are $\gamma$-a.i.'s. The converse is not quite true: the fixed point
set of a Hamiltonian $S^1$-action can have positive dimension.

Regarding Question \ref{quest:sa-ai}, it is also worth pointing out
that conjecturally the $\gamma$-diameter of the group of Hamiltonian
diffeomorphisms of a symplectically aspherical manifold $M$ is
infinite. In some instances, e.g., for surfaces, this is
obvious. Moreover, one might expect that $\gamma(\varphi^k)$ is
unbounded for generic (all?) maps $\varphi\neq \id$ or even that
$\gamma(\varphi^k)\to\infty$.

One interesting class of $\gamma$-a.i.'s, relevant for what follows,
is identified in \cite{GG:PR}. These are (Hamiltonian)
pseudo-rotations of $\CP^n$, i.e., Hamiltonian diffeomorphisms of
$\CP^n$ with minimal possible number of periodic points, equal to
$n+1$ by the Arnold conjecture, \cite{MS, Sa}. Among pseudo-rotations
are the Anosov--Katok pseudo-rotations from Example \ref{Ex:PR-2D} and
true rotations (i.e., isometries) of $\CP^n$ with finitely many fixed
points.  The following theorem is proved in a slightly different form
in \cite{GG:PR}:

\begin{Theorem}[$\gamma$-convergence, Thm.\ 5.1, \cite{GG:PR}]
\label{thm:gamma}
Let $\varphi$ be a pseudo-rotation of $\CP^n$. Then $\varphi$ is
$\gamma$-almost periodic. Furthermore, there exist a constant $C>0$
and a non-negative integer $d\leq n$, both depending only on
$\varphi$, such that for every $\eps$ in the range $(0,\,\pi)$, we
have
\[
\liminf_{k\to\infty} \frac{|\{\ell\leq k \mid 
\gamma(\varphi^\ell)<\eps\}|}{k}\geq C\eps^d.
\]
\end{Theorem}

The proof of the theorem is based on trading the behavior of
$\varphi^k$ with respect to the $\gamma$-norm for the dynamics of a
certain translation in the $d$-dimensional torus, which is of course
almost periodic; see Example \ref{ex:torus}.  We note that $\varphi$
cannot be $C^0$-almost periodic when it has a dense orbit -- the proof
of this fact is identical to the argument in Example
\ref{Ex:PR-2D}. However, $\varphi$ is a $C^0$-a.i.\ when it meets a
certain additional requirement generalizing the condition from that
example that the rotation number $\theta/\pi$ is exponentially
Liouville, \cite[Thm.\ 1.4]{GG:PR}. It is unknown if every Hamiltonian
pseudo-rotation of $\CP^n$ is a $C^0$-a.i.\ and if there are
$\gamma$-a.i.'s on $\CP^n$ which are not pseudo-rotations. (Note that
the $\gamma$-diameter of the universal covering of the group of
Hamiltonian diffeomorphisms of $\CP^n$ is bounded -- in fact equal to
$\pi$ -- in contrast with symplectically aspherical manifolds,
\cite{EP}.)  In any event, one can expect $\gamma$-a.i.'s to be
extremely non-generic for all $M$; for Hamiltonian diffeomorphisms
with finitely many periodic orbits this is proved in \cite{GG:generic}
for many symplectic manifolds including $\CP^n$.

We conclude this section by pointing out that Question
\ref{quest:interior} can be meaningfully restricted to Hamiltonian
transformations. For $C^0$-almost periodic compactly supported
Hamiltonian diffeomorphisms, just as in the general case of the
Hilbert--Smith conjecture, it boils down to the question whether the
action of the group of $p$-adic integers on $M$ can be generated by
such a transformation; cf., e.g., \cite{Pa, Ta}.

\section{Lagrangian Poincar\'e recurrence}
\label{sec:LPR}

In this section we will concentrate on an apparently different
question, which we prefer to state as a conjecture. The question is
connected to the $\gamma$-convergence theorem (Theorem
\ref{thm:gamma}), but this connection might be purely accidental.  As
in Section \ref{sec:AI-Ham}, let $\varphi$ be a compactly supported
Hamiltonian diffeomorphism of a symplectic manifold $M^{2n}$. The
following conjecture was put forth by the first author and
independently by Claude Viterbo around 2010.

\begin{Conjecture}[Lagrangian Poincar\'e Recurrence]
\label{conj:LPR}
For any closed Lagrangian submanifold $L\subset M$ there exists a
sequence of iterations $k_i\to\infty$ such that
\[
  \varphi^{k_i}(L)\cap L\neq \emptyset.
\]
Moreover, the density of the sequence $k_i$ is related to a symplectic
capacity of $L$.
\end{Conjecture}

We refer the reader to \cite[Sect.\ 5.1.2]{GG:PR} for a detailed
discussion of this conjecture. Here we only mention that the
conjecture is mainly interesting when $L$ is small (e.g., contained in
a small ball or more generally displaceable), and that the condition
that $\varphi$ is Hamiltonian is essential. In dimension two,
Conjecture \ref{conj:LPR} is an easy consequence of the standard
Poincar\'e recurrence. In general, very little is known about the
problem. It is not even clear if this conjecture is a dynamics
question or a packing problem: it is possible that there is an upper
bound on the number of disjoint Lagrangian submanifolds, Hamiltonian
diffeomorphic to each other and embedded into a compact domain in
$M$. However, at the time of writing, the only non-trivial result
along the lines of the conjecture is on the dynamics side. This is the
following theorem proved in \cite{GG:PR} and establishing a strong
form of the Lagrangian Poincar\'e recurrence for pseudo-rotations of
$\CP^n$.

\begin{Theorem}[Thm.\ 4.2, \cite{GG:PR}]
\label{thm:LPR}
Let $\varphi$ be a pseudo-rotation of $\CP^n$ and let $L\subset \CP^n$
be a closed Lagrangian submanifold, which admits a Riemannian metric
without contractible closed geodesics (e.g., a torus). Then
$\varphi^{k_i}(L)\cap L\neq \emptyset$ for some quasi-arithmetic
sequence $k_i\to\infty$. Furthermore, there exists a constant $C>0$
and a non-negative integer $d\leq n$, both depending only on $\varphi$
but not $L$, and a constant $a>0$ depending only on $L$ such that
\[
  \liminf_{k\to\infty} \frac{|\{\ell\leq k \mid \varphi^\ell(L)\cap
    L\neq \emptyset\}|}{k}\geq C\cdot a^d.
\]
\end{Theorem}

The constant $a$ is a certain homological capacity of $L$. Theorem
\ref{thm:LPR} is an easy consequence of Theorem \ref{thm:gamma} and
the standard fact that $\varphi(L)\cap L\neq\emptyset$ when
$\gamma(\varphi)<a$. In fact, the condition on $L$ in Theorem
\ref{thm:LPR} can be removed; see \cite[Thm.\ D and Rmk.\ 13]{KS}.
Note also that here and in Theorem \ref{thm:gamma} one could slightly
refine the quantitative result by replacing the lower bound on the
density of the sequence of ``returns'' by the upper bound $1/(C a^d)$ on
the step of the quasi-arithmetic sequence.

Conjecture \ref{conj:LPR} has deep applications to dynamics. For
instance, as has been pointed out by Viterbo, once established in
dimension four, it would imply (via a non-trivial argument) that the
group of area-preserving transformations of the closed disk is not
simple -- a well-known open question in two-dimensional dynamics.

\end{document}